\documentclass[12pt,a4paper]{article}
\usepackage{amstext}
\usepackage{fancyhdr}
\usepackage{amsfonts,graphicx,bezier, amssymb}
\usepackage{amsmath}
\usepackage{caption}
\newenvironment{prf}{\noindent{\bf{Proof:}}~~}{\hfill\rule{1ex}{1ex}\vskip1.5ex}
\newcommand{\Z}{\mathbb Z}

\newcommand{\N}{\mathbb N}

\newcommand{\beqa}{\begin{eqnarray}}
\newcommand{\enqa}{\end{eqnarray}}
\newcommand{\beq}{\begin{eqnarray*}}
\newcommand{\enq}{\end{eqnarray*}}

\newtheorem{qn}{Question}[section]
\newtheorem{rem}{Remark}[section]

\newtheorem{cor}{Corollary}[section]
\newtheorem{propn}{Proposition}[section]
\newtheorem{defn}{Definition}[section]
\newtheorem{exam}{{\bf Example}}[section]
\newtheorem{thm}{Theorem}[section]

\newtheorem{lem}{Lemma}[section]

\newcommand{\noi}{\noindent}

\begin{document}

\begin{center}\Large{\bf{Effect of nilpotency on semisimplicity and cohomology  of the $\Z$-module $\Z/(p_1^{k_1}\times \cdots \times p_n^{k_n})\Z$}}

\end{center}
\vspace*{0.3cm}
\begin{center}

David Ssevviiri\\

\vspace*{0.3cm}
Department of Mathematics\\
Makerere University, P.O BOX 7062, Kampala Uganda\\
Email: ssevviiri@cns.mak.ac.ug, ssevviirid@yahoo.com
\end{center}

\begin{abstract}
We show that existence of nonzero nilpotent elements in the $\Z$-module $\Z/(p_1^{k_1}\times \cdots \times p_n^{k_n})\Z$ inhibits the module from 
possessing good structural properties. In particular, it stops it from being semisimple and from  admitting certain good homological properties.
   \end{abstract}

{\bf Keywords}:    Semisimple modules; Cohomology of modules; Reduced modules; Nilpotent elements of modules

\vspace*{0.4cm}

{\bf MSC 2010} Mathematics Subject Classification:16D70,   16D60, 13D03

\section{Introduction}
\begin{paragraph}\noi 
 Nilpotent elements of a ring play a big role in studying properties of that ring. The collection of all  nilpotent
  elements of a commutative ring  is the intersection of all prime ideals of that ring. The sum of all nil ideals of a ring (ideals whose
   all elements are nilpotent) forms a radical called the K\"{o}the upper nil  radical.
   Nilpotent elements of a not necessarily commutative
    ring are used to define  a popular class of rings called 2-primal rings. A ring is {\it 2-primal} if the collection of all its
     nilpotent elements coincides with its prime radical. Commutative rings and reduced rings are 2-primal.
  \end{paragraph}
  
  \begin{paragraph}\noi     
     For modules, nilpotent elements as defined in \cite{SsevviiriNico} also reveal some properties about a module. In this
      paper,  we  show  that their existence, inhibits a module from possessing certain structural properties.
      In particular, a $\Z$-module $M=\Z/(p_1^{k_1}\times \cdots\times p_n^{k_n})\Z$ is semisimple if and only if it has no nonzero
      nilpotent elements and such a module exhibits ``good'' homological properties whenever it has no nonzero
      nilpotent elements.
\end{paragraph}

\begin{paragraph}\noi
All rings are unital and commutative;  the modules are unital left-modules over commutative rings. By drawing motivation
from how nilpotent 
 elements are defined in  rings, nilpotent elements of modules were  defined in \cite{SsevviiriNico}. This definition proved 
  instrumental in the papers:  \cite{Daudi} and \cite{SsevviiriNico}. In \cite{SsevviiriNico}, it was generalized
  to strongly nilpotent elements. Strongly nilpotent elements of a uniserial module coincide with the classical prime radical of that module,
   see \cite[Theorem 3.1]{SsevviiriNico} and \cite{Corrigendum}. This result generalizes that of Levitzki for rings, see
   \cite[Theorem 2.6]{McConnell}. In   \cite{Daudi}, it was used to characterize  non-nilpotent elements of the 
   $\Z$-module   $\Z/(p_1^{k_1}\times \cdots \times p_n^{k_n})\Z$. The direct limit of    
    non-nilpotent elements of this    module adjoined with  the zero element forms  an ideal  of 
    a direct sum of the ring of $p$-adic integers,         \cite[Theorem  1.2]{Daudi}. 
    
 \end{paragraph}

\section{Nilpotent elements and Reduced modules}
\begin{defn}\label{1}\rm\cite[Definition 2.1]{SsevviiriNico}
  An element $m$ of an $R$-module $M$ is {\it nilpotent}
  if $m=0$ or there exists $r\in R$ and  a positive integer  $k$ greater than 1 such that $r^km=0$ but
$rm\not=0$.
\end{defn}

\begin{exam}\label{e1}\rm
 Let $k$ be a field,  $A=k[x]/(x^2)$ an algebra and $M$ a 2-dimensional module defined over $A$. $A$ has basis $\{1, x\}$. Let the basis of 
$M$ be $\{u, v\}$ and the  action of  (basis of) $A$ on (the basis of) $M$ be given by 
\begin{center}
 
\begin{tabular}{c|cc}
  $A\setminus M$ & $u$ &  $ v$\cr\hline
   1 & $u$ & $v$\cr
   $x$ & $v$ & $0$. \cr
\end{tabular}
\end{center}

Any $m=pu+qv\in M$  with $0\not=p\in k$ such that $v\not=0$ is nilpotent. For if $0\not=a\in A$ is given by $a=rx$ for some $r\in A$, it
 follows that  $am= (rx)(pu+qv)= rpv\not=0$ but $a^2m=0$.
  \end{exam}

 \begin{exam}\rm\label{Ex2}
  
 In the $\Z$-module $M=\Z/(p_1^{k_1}\times \cdots \times p_n^{k_n})\Z$,
  all elements of $M$ which are not of the form $\{mp_1^{k_1-1}\cdots p_n^{k_n-1}\}_{m=1}^{(p_1p_2\cdots p_n)-1}$  are nilpotent.
\end{exam}

  \begin{paragraph}\noi 
   The number of non-nilpotent elements in the $\Z$-module $\Z/(p_1^{k_1}\times \cdots \times p_n^{k_n})\Z$ is $p_1p_2\cdots p_n - 1$. This number  is 
not dependent on the powers of the primes but rather   on the primes in question. The number $p_1p_2\cdots p_n - 1$ is an invariant for 
modules involving the same primes (they could have different multiplicities).
We say that the $\Z$-modules of the type above are in the same class if they have the same number of non-nilpotent
elements.
   \end{paragraph}

\begin{defn}\rm 
An $R$-module $M$ is {\it reduced}   if for all $r\in R$ and every $m\in M$, $r^2m=0$ implies $rm=0$. 
\end{defn}

\begin{paragraph}\noi
It is evident that a module is
 reduced if and only if it has no nonzero nilpotent elements.
   \end{paragraph}

\begin{exam}\rm
 If $G$  is a free abelian group, then when considered as a $\Z$-module, it has no nonzero nilpotent elements and hence
 it is reduced.
\end{exam}

\section{Effect on Semisimplicity}
\begin{lem}\label{LL1} A simple module is reduced.
\end{lem}

\begin{prf}
  An $R$-module $M$ is {\it prime} if   $rm=0$ for $r\in R$ and $m\in M$ implies either $m=0$ or $rM=0$. 
  We show that a simple module is prime. Let $rm=0$ for $m\in M$ and $r\in R$ where $M$ is a simple $R$-module.
  If $m=0$, $M$ is a prime module. Suppose $m\not=0$. Then $Rm\not=\{0\}$. So, $Rm=M$ since $M$ is simple. From $rm=0$ we get $rRm=0$ so that $rM=\{0\}$.
  Hence,   $M$ is prime. It is easy to see that any prime module  is reduced.
\end{prf}

\begin{paragraph}\noi 
A module is {\it semisimple} if it is a direct sum of simple modules.  
 \end{paragraph}
 
 \begin{propn}\label{P1}
  Any semisimple module is reduced.
 \end{propn}
 
 \begin{prf}
  Suppose $M$ is a semisimple $R$-module, $m\in M$ and $r\in R$ such that $r^2m=0$. Then $r^2(m_1, m_2, \cdots, m_i, \cdots)=0$ where
  $m=(m_1, m_2, \cdots, m_i, \cdots)\in \oplus_{i\in I}M_i=M$
   for some simple modules  $M_i$. So, $r^2m_i=0$ for each $i$. Since a simple module  is reduced
    (by Lemma \ref{LL1}), we have    $rm_i=0$ for each $i$. It follows that $rm=0$ and hence $M$ is reduced. 
 \end{prf}

\begin{rem}\rm\label{R1}  
 Proposition \ref{P1} provides us with a criterion for non-semisimplicity of   modules over commutative rings.  
   Any module with at least one nonzero nilpotent element cannot be semisimple. This fact also holds for rings. We know that  a ring $R$ is reduced 
  if and only if $_RR$ is a reduced module and, a ring $R$ is semisimple if and only if $_RR$ is semisimple. So, we have   
  $$ R ~\text{is semisimple}\Leftrightarrow~~ _RR ~\text{is semisimple}\Rightarrow~~ _RR ~\text{is reduced } \Leftrightarrow R ~\text{is reduced}.$$ 
  Thus, a ring    which contains at least one nonzero nilpotent element is not semisimple.
 \end{rem}

A reduced module need not be semisimple. The $\Z$-module $\Z$ is reduced but not semisimple.

\begin{thm}\label{T1}
 The following statements are equivalent for the $\Z$-module $M=\Z/(p_1^{k_1}\times \cdots \times p_n^{k_n})\Z$:
 \begin{enumerate}
  \item  $M$ is semisimple,
  \item $M$ is reduced,
  \item $k_1=\cdots=k_n=1$,
   \item $M$ is  a torsion module and if $p^2m=0$ for some prime number $p$ and $m\in M$, then $pm=0$.
 \end{enumerate}
\end{thm}

\begin{prf} $1\Leftrightarrow 3$. Since
\begin{equation}\label{equation1}
 \Z/(p_1^{k_1}\times \cdots \times p_n^{k_n})\Z\cong\Z/p_1^{k_1}\Z    \times \cdots \times \Z/p_n^{k_n}\Z, 
\end{equation}

$M$ is semisimple if and only if the $\Z$-modules   on the right hand side of the isomorphism (\ref{equation1}) are simple. The $\Z$-modules 
$\Z/p_1^{k_1}\Z, \Z/p_1^{k_1}\Z,  \cdots, \Z/p_n^{k_n}\Z$ are all simple if and only if $k_1=k_2=\cdots =k_n=1$.
It follows that $M$ is semisimple if and only if $k_1=k_2=\cdots =k_n=1$. \\

$2\Leftrightarrow 3$. The non-nilpotent elements of the 
$\Z$-module $M=\Z/(p_1^{k_1}\times \cdots \times p_n^{k_n})\Z$ are
 \begin{equation}
 \{mp_1^{k_1-1}\cdots p_n^{k_n-1}\}_{m=1}^{(p_1p_2\cdots p_n)-1}, 
 \end{equation}
 see \cite[p.700]{Daudi}. We observe that $k_1=k_2=\cdots =k_n=1$ if and only if
  the would be   non-nilpotent elements of $M$ are the integers $\{m\}_{m=1}^{(p_1p_2\cdots p_n)-1}$ which coincide with the nonzero elements of 
  $\Z/(p_1\times \cdots \times p_n)\Z$, i.e., 
   if and only if $M$ is reduced.
  $1\Leftrightarrow 4$ follows from \cite[Ex. 20.17(7), p. 173]{Wisbauer}.
\end{prf}

\begin{cor}\label{CC}
 The following statements are equivalent for the $\Z$-module $M=\Z/(p_1^{k_1}\times \cdots \times p_n^{k_n})$:
 \begin{enumerate}
  \item  $M$ is not semisimple,
  \item $M$ has at least one nonzero nilpotent element,
  \item at least one $k_i$  is greater than one for $i\in \{1, \cdots n\}$,
    \end{enumerate}

\end{cor}

\section{Effect on  Cohomology}

\begin{thm}\label{t}
 Let $\Z_{p^k}$ denote the $\Z$-module $\Z/p^k\Z$ where $p$  is a prime integer and $k\in  \N$. Let $f_{p^i}$ denote the $\Z$-endomorphism 
 of $\Z/p^k\Z$ given   by left multiplication of $p^i$, $i\in \N$. Then 
   the sequence $$ \Z_{p^k} \stackrel{f_{p^r}}{\rightarrow}\Z_{p^k} \stackrel{f_{p^{r+1}}}{\rightarrow}
  \cdots \rightarrow \Z_{p^k} \stackrel{f_{p^{k-1}}}{\rightarrow}\Z_{p^k} 
  \stackrel{f_{p^k}}{\rightarrow}\Z_{p^k}    \stackrel{f_{p^{k+1}}}{\rightarrow}\Z_{p^k} \rightarrow \cdots$$ where $r$  is the least positive integer
  such that $2r+1\geq k$  (i.e., $r\geq \frac{k-1}{2}$)   is a co-chain complex whose  $n$th cohomology group ($n\not=0$) is  given by 
  $$H^n(\Z_{p^k})=\begin{cases}
  \Z_p & k=1\\
  \frac{p^{(k-1-n)}\Z_{p^k}}{p^n\Z_{p^k}} & \frac{k}{2}\leq n \leq k-1, \text{for}~k~\text{even and;}~\frac{k-1}{2}\leq n \leq k-1,
  \text{for}~1< k~\text{odd} \\
  \Z_{p^k} & 1< k\leq n
                   
                  \end{cases}
$$

\end{thm}
 
 \begin{prf}
  Routine.
 \end{prf}

\begin{cor}
 Consider the hypotheses in Theorem \ref{t}. Then 
 \begin{enumerate}
  \item whenever $k=2n+1$ $\text{for}~k~\text{even and;}~\frac{k-1}{2}\leq n \leq k-1,
  \text{for}~1< k~\text{odd}$, we get $$H^n(\Z_{p^k})=\{0\};$$
  \item the sequence $\{ H^n(\Z_{p^k})\}_{n=1}^{\infty}$ for $k> 1$ is convergent and $$\lim_{ n\to\infty}    H^n(\Z_{p^k})= \Z_{p^k};$$
  \item the number of cohomology groups distinct from $\Z_{p^k}, k>1$ is  an even number given by  
  $$\begin{cases}
     \frac{k}{2} & k  ~\text{even} \\
     \frac{k+1}{2} & k ~\text{odd}.
    \end{cases}
  $$
  \item Computation of $H^n(\Z_{p_1^{k_1}\times \cdots \times  p_n^{k_n}})$ is easy since 
  $$\Z_{p_1^{k_1}\times \cdots \times  p_n^{k_n}}\cong  \Z_{p_1^{k_1}}\times \cdots \times  \Z_{p_n^{k_n}}$$ and 
  cohomology commutes with finite direct sum.
 \end{enumerate}

\end{cor}

\begin{rem}\rm
 Whereas the number of non-nilpotent elements of $\Z_{p^k}$ is dependent on $p$ and independent of $k$, the number of cohomology groups depends on $k$ 
 and is independent of $p$.
\end{rem}

Take a $\Z$-module $M=\Z/(2^3\times 3^2\times 5^3)\Z$. $M\cong \Z/2^2\Z \times \Z/3^2\Z \times \Z/5^3\Z$.
So, $H^1(M)= 0 \times \frac{\Z_9}{3\Z_9} \times 0$, 
$H^2(M)= \frac{\Z_8}{4\Z_8} \times \Z_9 \times \frac{\Z_{125}}{25\Z_{125}}$, $H^3(M)= \Z_8 \times  \Z_9  \times \Z_{125}$, and 
$H^m(M)=H^3(M)$ for all $m\geq 3$.\\

In general,  existence of nonzero nilpotent elements of    a $\Z$-module $M=\Z/(p_1^{k_1}\times \cdots \times p_n^{k_n})\Z$   inhibits
the sequence $\{H^n(M)\}_{n=1}^{\infty}$ from being a constant sequence and  the $k$s (i.e., powers to the primes $p_1, \cdots, p_n$) determine the 
rate at which $\{H^n(M)\}_{n=1}^{\infty}$ becomes constant. In particular, the number of steps required for the sequence 
$\{H^n(M)\}_{n=1}^{\infty}$ to become stable is $\text{Maximum}\{k_1, \cdots, k_n\}$. Furthermore, under the sequence $\{ H^n(-)\}_{n=1}^{\infty}$,
$\Z$-modules $M=\Z/(p_1^{k_1}\times \cdots \times p_n^{k_n})\Z$  with at least one $k_i> 1$  stabilize  to the same module 
if and only if they belong to the same class (i.e., they posses  the same number of non-nilpotent elements). We are then led to the following theorem.

\begin{thm} For a $\Z$-module   $M=\Z/(p_1^{k_1}\times \cdots \times p_n^{k_n})\Z$, the following statements are equivalent:
\begin{enumerate}
 \item $M$  is reduced,
 \item $\{H^n(M)\}_{n=1}^{\infty}$  is a constant sequence given by $\Z/p_1\Z \times \Z/p_2\Z \times \cdots \times \Z/p_n\Z$,
 \item the $\Z$-module $\Z/p_1\Z \times \Z/p_2\Z \times \cdots \times \Z/p_n\Z$ is invariant under the cohomology $H^n(-)$ for all $n\in \Z^+$,
 \item $k_1= k_2, \cdots =k_n =1$.
\end{enumerate}

\end{thm}

 If 
  $M=\Z/(p_1^{k_1}\times \cdots \times p_n^{k_n})\Z$ where at least one $k_i>1$, then $M$ is not invariant under $H^n(-)$ for all $n\in \Z^+$.
  The cause of its failure to be invariant is the existence of  nonzero nilpotent elements.

\section{Getting rid of the bad elements}

\begin{paragraph}\noi 
 One of the central  ideas of radical theory in associative rings, is to collect all   ``bad'' elements of a ring $R$  in some ideal $\gamma (R)$
 called the radical ideal of $R$  such that   the ring $R/\gamma(R)$ is ``good'', i.e.,    admits a ``good'' structural theorem, see \cite[p. 21]{Gardner}.
 For instance, if $\mathcal{N}(R)$ is the collection of all nilpotent elements of a ring $R$ and nilpotent elements are considered ``bad''
 elements, then $R/\mathcal{N}(R)$ is ``good'', i.e., contains no bad elements since  $\mathcal{N}(R/\mathcal{N}(R))=\{0\}$.
 \end{paragraph}
 
 \begin{paragraph}\noi
 We observe that, the above strategy does not work for a $\Z$-module $M=\Z/(p_1^{k_1}\times \cdots \times p_n^{k_n})\Z$ where at least one $k_i>1$.
 This is because: 1) whereas for rings $\mathcal{N}(R)$ is an ideal of $R$, for modules, $\mathcal{N}(M)$ need not be a submodule of $M$;
 2) if we take  the smallest submodule of $M$  containing  $\mathcal{N}(M)$   denoted by $\langle  \mathcal{N}(M)\rangle$, it follows
 that $\langle  \mathcal{N}(M)\rangle=M$ so that on factoring, nothing is left.
  \end{paragraph}
  
  \begin{paragraph}\noi 
So, a submodule $N$ of $M$ is chosen so that  $M/N$  admits   
  our desired structural properties, (i.e., $M/N$ is semisimple and is
   invariant under $H^n(-)$ for all $n\in \Z^+$) and this factor module is the largest possible.
 \end{paragraph}

  \begin{paragraph}\noi 
   
 For instance, if $M=\Z/({p_1}^2 \times {p_2}^4\times p_3)\Z$, then $N= \Z/({p_1}\times {p_2}^3)\Z$ so that
  $$M/N\cong \Z/(p_1\times p_2\times p_3)\Z.$$ In general, for $M=\Z/(p_1^{k_1}\times \cdots \times p_n^{k_n})\Z$, 
  $N=\Z/(p_1^{k_1-1}\times \cdots \times p_n^{k_n-1})\Z$.

  \end{paragraph}

   \begin{propn}
    Let   $M=\Z/(p_1^{k_1}\times \cdots \times p_n^{k_n})\Z$ be a $\Z$-module and $N$ be the set of all non-nilpotent elements
    of $M$ together with the zero element. Then $$M/N\cong \Z/(p_1^{k_1-1}\times \cdots \times p_n^{k_n-1})\Z.$$
   \end{propn}

   \begin{prf}
    Let $N=\{0\}\cup \{\text{non-nilpotent elements of}~M\}$. Then $$N=\{0\}\cup \{mp_1^{k-1}\cdots p_n^{k_n-1}\}_{m=1}^{p_1p_2\cdots p_n -1}$$ So,
    $$N=p_1^{k-1}\cdots p_n^{k_n-1}\{0, 1, 2, \cdots, p_1p_2\cdots p_n -1\}.$$ By \cite[Proposition 1.4]{Daudi}, $N$  is isomorphic
    to $\Z/(p_1, \cdots, p_n)\Z$. Hence $$M/N\cong \Z/(p_1^{k_1-1}\times \cdots \times p_n^{k_n-1})\Z.$$
   \end{prf}

   \begin{paragraph}\noi 
   This shows that successive factoring out of  non-nilpotent elements  adjoined with the zero element eventually leaves no non-nilpotent elements.
     \end{paragraph}

\begin{qn}
Is there a   general setting  in which  the above results hold?
\end{qn}

\section*{Conclusion}

Nilpotent elements control the structure of the $\Z$-module $\Z/(p_1^{k_1}\times \cdots \times p_n^{k_n})\Z$.

\end{document}